\documentstyle[12pt]{article}
\newtheorem{Th}{Theorem}
\newtheorem{Prop}[Th]{Proposition}
\newtheorem{Def}[Th]{Definition}
\newtheorem{Rem}[Th]{\it Remark}
\newtheorem{Cor}[Th]{Corollary}

\newcommand{\Proof}{\noindent {\it Proof. }}
\newcommand{\be}{\begin{eqnarray*}}
\newcommand{\ee}{\end{eqnarray*}}

\newcommand{\pkXY}{{\cal P}(^k\!X;Y)}

\newcommand{\pk}{{\cal P}(^k\!}

\newcommand{\pwbkXY}{{\cal P}_{wb}(^k\!X;Y)}

\newcommand{\pwbk}{{\cal P}_{wb}(^k\!}
\newcommand{\finesp}{\hspace*{\fill} $\Box$\vspace{0.5\baselineskip}}
\newcommand{\fin}{\hspace*{\fill} $\Box$}

\newcommand{\LkXY}{{\cal L}(^k\!X;Y)}
\newcommand{\LwbkXY}{{\cal L}_{wb}(^k\!X;Y)}
\newcommand{\LXkY}{{\cal L}(X_1,\ldots ,X_k;Y)}
\newcommand{\LZkY}{{\cal L}(Z_1,\ldots ,Z_k;Y)}
\newcommand{\LwbXkY}{{\cal L}_{wb}(X_1,\ldots ,X_k;Y)}
\newcommand{\LwbZkY}{{\cal L}_{wb}(Z_1,\ldots ,Z_k;Y)}

\newcommand{\HXY}{{\cal H}(X;Y)}
\newcommand{\HwbuXY}{{\cal H}_{wbu}(X;Y)}

\newcommand{\N}{\bf N}

\newcommand{\lra}{\longrightarrow}
\newcommand{\Ra}{\Rightarrow}

\newcommand{\ra}{\rightarrow}
\newcommand{\eps}{\epsilon}
\newcommand{\noin}{\noindent}
\newcommand{\espv}{\vspace{.5\baselineskip}}
\newcommand{\espvv}{\vspace{\baselineskip}}

\begin{document}

\title{Factorization of weakly continuous \\
      holomorphic mappings}

\author{Manuel Gonz\'alez\thanks{Supported in part by DGICYT Grant PB
      91--0307 (Spain)}
       \and
      Joaqu\'\i n  M. Guti\'errez\thanks{Supported in part by DGICYT
      Grant PB 90--0044 (Spain)}}

\date{\hspace*{\fill}}
\maketitle

\begin{abstract}
We prove a basic property of continuous multilinear mappings
between topological vector spaces, from which we derive an easy proof
of the fact that a multilinear mapping (and a polynomial) between
topological vector spaces is weakly continuous on weakly bounded sets if and
only if it is weakly {\it uniformly\/} continuous on weakly bounded sets.
This result was obtained in 1983 by Aron, Herv\'es and Valdivia
for polynomials between Banach spaces, and it also holds if
the weak topology is replaced by a coarser one.
However, we show that it need not be true for a stronger topology,
thus answering a question raised by Aron.
As an application of the first result,
we prove that a holomorphic mapping $f$ between complex Banach
spaces is weakly uniformly continuous on bounded subsets if and only if
it admits a factorization of the form $f=g\circ S$, where $S$
is a compact operator and $g$ a holomorphic mapping.
\end{abstract}

\vspace{\fill}

\newpage

Our aim is to give characterizations of polynomials and holomorphic
mappings on Banach spaces, which are weakly uniformly
continuous on bounded sets.
The polynomials with this property have been studied
by many authors: see, for instance,
\cite{AGL,AHV,AP}.
A reason for their interest is that
they are uniform limits of finite type
polynomials (assuming the approximation property on the dual space)
\cite[Proposition~2.7]{AP}.
In the case of locally convex spaces, these classes of polynomials
have been analysed in several places
(see, e.g., \cite{Mo}).
The holomorphic mappings with weakly uniformly continuous restrictions
to bounded sets have also been considered in various papers
\cite{AHV,Dic0,Mo}.

The paper is organized in three sections.
In the first one we prove a basic result (Theorem~\ref{klinear}
below) on continuity of multilinear mappings between topological vector
spaces, roughly showing that a multilinear mapping is continuous on
certain subsets if and only if the Cauchy nets contained in these
subsets are mapped into Cauchy nets.

As an easy consequence, we show that if a multilinear mapping between
topological vector spaces is weakly continuous on weakly bounded sets,
then it is also {\it uniformly\/} continuous on them.
This is also true if the weak topology is replaced by a coarser
one, and the proof uses the fact that weakly bounded sets are
weakly precompact.

This result extends and simplifies a well known theorem proved by
Aron, Herv\'es and Valdivia \cite[Theorem~2.9]{AHV}, for polynomials
between Banach spaces.
In fact, to prove Theorem~\ref{klinear}, we just
refine what they did in \cite{AHV}.
This first Section is quite helpful to understand what it means
for a multilinear mapping to be continuous or uniformly continuous
on certain classes of subsets.

In the second Section, we answer negatively a question of Richard
Aron's which was open for a number of years.
Namely, he asked if given any vector topology $\tau$ on a Banach
space $X$, it is true that
$\tau$-continuity on bounded sets is always equivalent
to uniform $\tau$-continuity on bounded sets, for a polynomial on $X$.
This was known to be the case when $\tau$ is the norm topology, the
weak topology, or the weak-star topology on a dual space.
We show
that if the unit ball is not $\tau$-precompact, then uniform
$\tau$-continuity does not necessarily follow from $\tau$-continuity.
To this end, we consider a locally convex topology, called $ckw$, defined
as the finest locally convex topology on a Banach space
having the same convergent sequences
as the weak topology.
This topology has been studied in various contexts \cite{GGkw,Po,We}.
Since the unit ball of $L_1[0,1]$ is not $ckw$-precompact,
and we possessed a nice decription of the $ckw$-topology on this
space, we were able
to construct a polynomial on $L_1[0,1]$ which is $ckw$-continuous on
bounded sets, but is not uniformly $ckw$-continuous on bounded sets.

Before giving the contents of Section \ref{factres},
we recall two easy properties of weakly uniformly continuous mappings
between Banach spaces:
first, every mapping whose restrictions to bounded sets
are weakly uniformly
continuous, takes bounded sets into
relatively compact sets \cite[Lemma~2.2]{AP},
and second, a linear bounded operator is compact
if and only if it is weakly (uniformly) continuous on bounded sets
\cite[Proposition~2.5]{AP}.

It is then clear that if a holomorphic mapping
$f$ between complex Banach spaces
admits a factorization of the form $f=g\circ S$, where $S$ is a compact
operator, and $g$ a holomorphic mapping, then $f$ is weakly uniformly
continuous on bounded sets.
In Section~\ref{factres}, we apply the main result of the
first part (Theorem~\ref{klinear}) to
show that these easy examples are the {\it only\/} ones,
i.e., every holomorphic mapping whose restrictions to bounded sets are
weakly uniformly continuous admits a factorization as above.
We take advantage of work by Braunsz and Junek: namely,
some ideas of Theorem~\ref{facklin}
below are contained in \cite[Proposition~2.14]{BJ}.

Factorizations of holomorphic mappings have already been considered,
but here the factors stand in inverse order.
Thus, it is proved in \cite{AS} that a holomorphic mapping $f$ is
compact (i.e., takes a neighbourhood of each point into a relatively
compact set) if and only if it admits a factorization of the form
$f=S\circ g$, where $g$ is a holomorphic mapping, and $S$ a compact
operator.
A similar result is proved in \cite{Ry} for
weakly compact holomorphic
mappings.\espv

If $X$, $Y$ are topological vector spaces,
to each $k$-homogeneous
polynomial $P$ from $X$ into $Y$,
it is associated a unique symmetric $k$-linear mapping $\hat{P}:
X\times\stackrel{(k)}{\ldots}\times X\ra Y$,
given by the polarization formula \cite[Theorem~1.10]{Mu}, so that
$\hat{P}(x,\ldots,x) = P(x)$ for all $x\in X$.
We refer to \cite{Din,Mu} for the general theory of polynomials and
holomorphic mappings on infinite dimensional spaces.

The set of natural numbers is denoted by $\N$.
A net in a topological vector space is said to be {\it (weakly) null\/}
if it is (weakly) convergent to zero.
If $X$ is a Banach space, then $B_X$ stands for its closed unit ball,
and $X^*$ for its dual.
By an {\it operator\/} we mean a linear continuous mapping.

\section{Uniformly continuous polynomials}

      In this part, we give a basic property of multilinear mappings,
and derive an easy proof of the equivalence between weak continuity
and weak uniform continuity on weakly bounded sets.

\begin{Def}
{\rm A family ${\cal B}$ of subsets of a vector space
$X$ is said to be a
{\it bornology\/} if it satisfies the following conditions:

(a) ${\cal B}$ covers $X$;

(b) $A\in {\cal B}, \; D\subset A \Ra D\in {\cal B}$;

(c) $A,D\in {\cal B} \Ra A\cup D\in{\cal B}$;

(d) for every $A\in {\cal B}$, and every scalar $\lambda$, we have
      $\lambda A\in{\cal B}$;

(e) $A,D\in{\cal B} \Ra A+D\in {\cal B}$.
}
\end{Def}

If $X$ is a topological vector space, examples of bornologies on $X$ are
the family of all subsets, the family of all (weakly) bounded subsets,
the (weakly) compact sets, etc.
Given a bornology ${\cal B}$ on $X$,
we say that $A\subseteq X$ is a ${\cal B}$-{\it
set\/}
if $A\in {\cal B}$. We say that a net
$(x_\alpha)_{\alpha\in\Gamma}\subset X$ is a
${\cal B}$-{\it net\/} if $\{ x_\alpha:\alpha\in\Gamma\} \in {\cal B}$.

\begin{Def}
\label{bwcont}
{\rm
Given topological vector spaces $X_1,\ldots,X_k,Y$, a $k$-linear
(not necessarily continuous) mapping $A:X_1\times\cdots\times X_k\ra Y$,
and a bornology ${\cal B}_j$ on
$X_j$, for each $1\leq j\leq k$, we say that $A$ is {\it continuous on
${\cal B}_j$-sets\/} if for each $B_j\in{\cal B}_j$, $x^j\in B_j$
$(1\leq j\leq k)$, and each zero neighbourhood $V$ in $Y$, there are
zero neighbourhoods $U_j$ in $X_j$ so that
$$
A(y^1,\ldots,y^k) - A(x^1,\ldots,x^k) \in V
$$
whenever $y^j\in B_j$
satisfy $y^j-x^j\in U_j$, for $1\leq j\leq k$.}
\end{Def}

Clearly, $A$ is continuous on ${\cal B}_j$-sets if and only if, given a
convergent ${\cal B}_j$-net $x^j_\alpha\ra x^j$
$(\alpha\in\Gamma)$ in $X_j$, for each $1\leq j\leq k$, then the net
$\left( A\left(
x^1_\alpha,\ldots,x^k_\alpha\right)\right)_{\alpha\in\Gamma}$
converges to $A(x^1,\ldots,x^k)$ in $Y$.

\begin{Th}
\label{klinear}
Let $X_1,\ldots,X_k,Y$ be topological vector spaces,
and $A:X_1\times\cdots\times X_k\ra Y$
a $k$-linear (not necessarily continuous) mapping.
Let ${\cal B}_j$ be a bornology on $X_j$ $(1\leq j\leq k)$.
Then the following assertions are equivalent:

{\rm (a)} $A$ is continuous on ${\cal B}_j$-sets;

{\rm (b)} given Cauchy ${\cal B}_j$-nets $(x_\alpha^j)_{\alpha\in
\Gamma}\subset
X_j$
$(1\leq j\leq k)$ such that at least one of them is null, then
the net $\left( A\left( x^1_\alpha,\ldots,x^k_\alpha\right)\right)
_{\alpha\in\Gamma}$ converges to zero in $Y$;

{\rm (c)}
given a Cauchy ${\cal B}_j$-net $\left( x^j_\alpha\right)_
{\alpha\in\Gamma}$ in $X_j$, for each
$1\leq j\leq k$, then the net
$\left( A\left( x^1_\alpha,\ldots,x^k_\alpha\right)
\right)_{\alpha\in\Gamma}$ is Cauchy in $Y$.
\end{Th}

\Proof
(a) $\Ra$ (b).
If $k=1$, there is nothing to prove.
Assume the result is true for all $(k-1)$-linear mappings
and fails for the $k$-linear mapping $A$.
Then we can find Cauchy ${\cal B}_j$-nets $\left( x^j_\alpha\right)
_{\alpha\in\Gamma}\subset X_j$ $(1\leq j\leq k)$, at least one of which
is null (to simplify notation,
assume $x^1_\alpha\ra 0$), and
a zero neighbourhood $V_1$ in $Y$ such that
$$
A\left( x^1_\alpha,\ldots,x^k_\alpha\right) \notin V_1 \hspace{2em}
(\alpha\in \Gamma) \, .
$$
Let $V_2$ be a zero neighbourhood such that $V_2+V_2\subseteq V_1$.
For each fixed $\alpha\in\Gamma$, the mapping $Ax^k_\alpha$ given by
$$
Ax^k_\alpha (z^1,\ldots,z^{k-1}):= A(z^1,\ldots,z^{k-1},x^k_\alpha )
\hspace{1.5em} \left( z^j\in X_j; 1\leq j\leq k-1\right)
$$
is $(k-1)$-linear and takes convergent ${\cal B}_j$-nets into convergent nets.
By induction, there is $\kappa (\alpha)\in\Gamma$ so that
$$
A\left( x^1_\beta,\ldots,x^{k-1}_\beta,x^k_\alpha\right) =
Ax^k_\alpha \left( x^1_\beta,\ldots,x^{k-1}_\beta\right) \in V_2
\hspace{1.5em} (\beta \geq \kappa (\alpha)) \, .
$$
For every $\alpha\in\Gamma$, we have
\be
\lefteqn{A\left( x^1_{\kappa (\alpha)},\ldots,x^{k-1}_{\kappa (\alpha)},
x^k_{\kappa (\alpha)}-x^k_\alpha\right) =} \\
& & A\left( x^1_{\kappa (\alpha)},\ldots,x^k_{\kappa (\alpha)}\right) -
A\left(x^1_{\kappa (\alpha)},\ldots,x^{k-1}_{\kappa (\alpha)},x^k_\alpha
\right) \notin V_2\, .
\ee
Consider the vectors:
$$
y^j_\alpha: = \left\{ \begin{array}{ll}
                        x^j_{\kappa (\alpha)} & \mbox{if $1\leq j\leq k-1$} \\
                        x^k_{\kappa (\alpha)}-x^k_\alpha & \mbox{if $j=k$}
                      \end{array} \right.
$$
We can assume $\kappa (\alpha) \geq \alpha$.
This condition assures us that the ${\cal B}_j$-nets
$\left( y^j_\alpha\right)_{\alpha\in\Gamma} \subset X_j$
$(1\leq j\leq k)$ are Cauchy,
and at least two of them are null.
By repeating the process, we obtain ${\cal B}_j$-nets
$\left( z^j_\alpha\right)_{\alpha\in\Gamma}
\subset X_j$ $(1\leq j\leq k)$, all of them null, and a zero neighbourhood
$V_k$ in $Y$, so that
$$
A\left( z^1_\alpha,\ldots,z^k_\alpha\right) \notin V_k \hspace{2em}
(\alpha\in\Gamma)\, .
$$
This contradicts our assumption (a).

(b) $\Ra$ (c).
Let $\left( x^j_\alpha\right)_{\alpha\in\Gamma}
\subset X_j$ $(1\leq j\leq k)$ be Cauchy
${\cal B}_j$-nets. We have
\be
\lefteqn{A\left( x^1_\alpha,\ldots,x^k_\alpha\right) -
              A\left( x^1_\beta,\ldots,x^k_\beta\right) = } \\
&  & A\left( x^1_\alpha -x^1_\beta,x^2_\alpha,\ldots,x^k_\alpha\right)
      +A\left( x^1_\beta, x^2_\alpha
              -x^2_\beta,x^3_\alpha,\ldots,x^k_\alpha\right) + \cdots \\
&  & \mbox{} + A\left( x^1_\beta, \ldots, x^{k-1}_\beta,
              x^k_\alpha - x^k_\beta\right)\, .
\ee
In each of the above terms, letting $(\alpha,\beta)\in\Gamma\times\Gamma$
increase, all the nets
are Cauchy ${\cal B}_j$-nets, and at least one of them is null. Then
$$
\lim_{\alpha,\beta} \left[ A\left( x^1_\alpha,\ldots,x^k_\alpha\right) -
A\left( x^1_\beta,\ldots,x^k_\beta\right) \right]= 0\, ,
$$
and therefore
$\left( A\left( x^1_\alpha,\ldots,x^k_\alpha\right)\right) _\alpha$ is
a Cauchy net.

(c) $\Ra$ (a).
Suppose $A$ is not continuous on ${\cal B}_j$-sets.
Then we can find sets $B_j\in{\cal B}_j$, points $x^j\in B_j$, and a
zero neighbourhood $V$ in $Y$ so that for every zero neighbourhood
$U_j$ in $X_j$ $(1\leq j\leq k)$
there is $y_{U_j}\in B_j$, with $y_{U_j}-x^j\in U_j$
but
$$
A\left( y_{U_1},\ldots,y_{U_k}\right) - A(x^1,\ldots,x^k) \notin V\, .
$$
Let ${\cal U}_j$ be the family of all zero neighbourhoods in $X_j$, and
${\cal U}:= {\cal U}_1\times\cdots\times {\cal U}_k\times\N$, ordered
in the natural way.
For each $U=(U_1,\ldots,U_k,n)\in{\cal U}$, and $j\in\{ 1,\ldots,k\}$, let
$$
z^j_U:= \left\{ \begin{array}{cl}
              y_{U_j} & \mbox{ if $n$ is even} \\
              x^j     & \mbox{ if $n$ is odd.}
              \end{array} \right.
$$
Then $\left( z^j_U\right)_{U\in{\cal U}}$ is a Cauchy ${\cal B}_j$-net.
However, the net
$\left( A\left( z^1_U,\ldots,z^k_U\right)\right)_{U\in{\cal U}}$ is
not Cauchy.
\finesp

If, for each $j$, we take as
${\cal B}_j$ the bornology of all subsets of $X_j$,
the assertion (a) in the last Theorem
simply states that $A$ is continuous.

\begin{Def}
{\rm
Let $f:X\ra Y$ be a mapping between topological vector spaces,
and ${\cal B}$ a bornology on $X$.
We say that $f$ is {\it uniformly continuous on ${\cal B}$-sets\/} if
for every zero neighbourhood $V$ in $Y$, and every $B\in{\cal B}$
there is a zero neighbourhood $U$ in $X$ such that we have
$f(x) - f(y)\in V$ whenever $x,y\in B$ satisfy $x-y\in U$.
This definition may be adapted to multilinear mappings in an obvious way.}
\end{Def}

Recall that a subset $B$ of a topological vector space $X$ is {\it
precompact\/} if for every zero neighbourhood $U$ in $X$ there is a
finite set $M\subseteq B$ such that $B\subseteq M+U$.
It is well known that $B$
is precompact if and only if every net in $B$ has a Cauchy subnet
\cite[Theorem~6.32]{Ke}.
We shall use this fact in the proof of the following result,
which relates uniform continuity with
the properties considered in Theorem~\ref{klinear}.

\begin{Th}
\label{equivAHV}
Let $X$, $Y$ be topological vector spaces, and
${\cal B}$ a bornology of precompact sets in $X$.
Then a mapping $f:X\ra Y$ is uniformly continuous on ${\cal B}$-sets if
and only if it takes Cauchy ${\cal B}$-nets into Cauchy nets.
\end{Th}

\Proof
Let $f$ be uniformly continuous on ${\cal B}$-sets,
$(x_\alpha)\subset X$ a
Cauchy ${\cal B}$-net, and $V$ a zero neighbourhood in $Y$.
There is a zero neighbourhood $U\subset X$ so that whenever
$x_\alpha-x_\beta\in U$, then we have $f(x_\alpha)-f(x_\beta) \in V$.
Now, since $(x_\alpha)$ is Cauchy, there is $\alpha_0$ such that
$$
x_\alpha - x_\beta \in U \hspace{2em} (\alpha ,\beta \geq \alpha_0) \, ,
$$
and hence $\left( f(x_\alpha)\right)$ is Cauchy.

Conversely, assume $f$ is not uniformly continuous on ${\cal B}$-sets.
Then we can find $B\in {\cal B}$, and a zero neighbourhood $V\subset Y$
so that for every
zero neighbourhood $U\subset X$, there are $x_U,y_U\in B$, with
$x_U-y_U\in U$ and $f(x_U)-f(y_U)\notin V$.

Let ${\cal U}$ be the family of all zero neighbourhoods in $X$.
Since every ${\cal B}$-net has a Cauchy subnet, we can assume that the
nets $(x_U)_{U\in{\cal U}}$, $(y_U)_{U\in{\cal U}}$ are Cauchy.
Consider the set ${\cal W}:={\cal U}\times
\N$, ordered in the natural way.
To each $W=(U,i) \in {\cal W}$ we associate
$$
z_W:= \left\{ \begin{array}{ll}
              x_U & \mbox{if $i$ is even}\\
              y_U & \mbox{if $i$ is odd.}
              \end{array} \right.
$$
Then the ${\cal B}$-net $(z_W)$ is Cauchy. However, $\left( f(z_W)\right)$
is not Cauchy.
\finesp

Clearly, Theorem \ref{equivAHV} is also valid for multilinear mappings,
with obvious modifications.

\begin{Rem}
{\rm
If $f:X\ra Y$ satisfies the hypotheses of Theorem~\ref{equivAHV},
then $f(B)$ is precompact for each $B\in{\cal B}$.
The converse is not true.
Indeed, by modifying an example given in \cite[p.~82]{Ll},
we now construct a real valued function $f$ on a Banach space,
with the following conditions:

(a) $f$ is weakly continuous on bounded sets (it will even be weakly
continuous on the whole space);

(b) $f$ takes bounded sets into precompact sets, i.e., $f$ is bounded
on bounded sets;

(c) $f$ is not weakly uniformly continuous on bounded sets.

Let $X$ be a separable, nonreflexive Banach space. By James' theorem,
we can find $\phi\in X^*$, with $\| \phi\|=1$,
which does not attain its
norm on $B_X$. Let
$$
g(x): = \frac{1}{\phi (x)-1} \hspace{2em} (x\in B_X)\, .
$$
Since $X$ is normal for the weak topology, and $g$ is weakly continuous,
$g$ admits an extension $\tilde{g}$ to $X$ which is weakly continuous.
Since $\tilde{g}$ is unbounded on $B_X$, it is not weakly uniformly
continuous on $B_X$.
Therefore, there is $\delta >0$ so that for each convex, weak zero
neighbourhood $U$ in $X$ we can find $x,y\in B_X$ with $x-y\in U$
and $|\tilde{g}(x) - \tilde{g}(y)| >\delta$.
The segment $[x,y]$ is clearly contained in $B_X\cap (y+U)$.
Choose $\lambda >\pi\delta^{-1}$, and define
$$
f(x):= \sin \left( \lambda\tilde{g}(x)\right) \hspace{2em} (x\in X)\, .
$$
Clearly, $f$ is weakly continuous and bounded.
However, it is not weakly uniformly continuous on $B_X$ since we can
find $z\in [x,y]$ so that $|f(z)-f(y)| \geq 1$.
}
\end{Rem}

\begin{Cor}
\label{wucmm}
Let $X_1,\ldots,X_k,Y$ be topological vector spaces.
Let $\tau_j$ be a vector topology on $X_j$, coarser than
or equal to the weak topology, and ${\cal B}_j$ a bornology of\/
{\em weakly bounded}
sets on $X_j$ $(1\leq j\leq k)$.
Then a $k$-linear mapping from $X_1\times\cdots\times X_k$ into $Y$
is $\tau_j$-continuous on ${\cal B}_j$-sets
if and only if it is uniformly $\tau_j$-continuous on ${\cal B}_j$-sets.
\end{Cor}

\Proof
Since the weakly bounded sets coincide with the weakly precompact sets
\cite[Corollary~8.1.6]{J},
every ${\cal B}_j$-set is $\tau_j$-precompact.
Therefore, it is enough to apply Theorem~\ref{klinear} and
Theorem~\ref{equivAHV}.
\finesp

      Using the polarization formula, it is clear that a polynomial
$P$ between topological vector spaces
takes convergent ${\cal B}$-nets into convergent
nets if and only if so does $\hat{P}$, and that $P$ takes Cauchy
${\cal B}$-nets into Cauchy nets if and only if so does $\hat{P}$.
Therefore, we obtain:

\begin{Cor}
\label{wucp}
Let $X$, $Y$ be topological vector spaces.
Let $\tau$ be a vector topology on $X$ coarser than or equal to the weak
topology, and ${\cal B}$ a bornology on $X$ consisting of\/
{\em weakly bounded} sets.
Then a homogeneous polynomial from $X$ into $Y$ is
$\tau$-continuous on ${\cal B}$-sets if and only if it is
uniformly $\tau$-continuous on ${\cal B}$-sets.
\end{Cor}

The particular case when
${\cal B}$ is the bornology of bounded sets in a Banach space, and
$\tau$ the weak topology was proved in \cite[Theorem~2.9]{AHV}.
The result for
${\cal B}$ being the bornology of Rosenthal sets, or Dunford-Pettis sets
in a Banach space, and $\tau$ the weak topology,
was obtained in \cite[Proposition~3.6 and the comment after it]{GGdw}.

\section{A counterexample}

In this Section, we show that the r\^ole of precompactness in
Corollaries~\ref{wucmm} and~\ref{wucp} is essential.
We shall consider a topology on $L_1[0,1]$,
compatible with the dual pairing $\langle L_1[0,1],L_\infty [0,1]
\rangle$,
for which the
unit ball is not precompact, and give an example of a polynomial
not satisfying the conclusion of Corollary~\ref{wucp} for this topology.

If $X$ and $Y$ are Banach spaces,
the space of $k$-homogeneous (continuous)
polynomials from $X$ into $Y$ is denoted by $\pkXY$, and that of $k$-linear
(continuous) mappings from $X^k=X\times\stackrel{(k)}{\ldots}\times X$
into $Y$, by $\LkXY$.
If $Y$ is omitted, it is understood to be the scalar field.

    The $ckw$ topology \cite{GGkw} on a Banach space $X$ is the finest
locally convex
topology having the same convergent sequences as the weak topology.
On an infinite dimensional space, it is strictly finer than the
weak topology, and on a Banach space without the Schur property,
it is strictly coarser than the norm topology.
It is therefore a topology compatible with the pairing
$\langle X,X^*\rangle$, and so its
bounded sets are the
norm bounded sets.

    Given a Banach space $X$,
    a subset $K$ of its dual $X^*$ is an ($L$)-{\it set\/}
if, for every weakly null sequence $(x_n)\subset X$, we have
$$
\lim_n \sup_{\phi\in K} |\langle x_n,\phi\rangle | = 0\, .
$$

    The $ckw$ topology on a Banach space turns out to be the topology
of uniform convergence on ($L$)-subsets of the dual
\cite[Theorem~3.1]{GGkw}. A subset is $ckw$-precompact if and only
if each sequence in it has a weak Cauchy subsequence
\cite[Theorem~4.4]{GGkw}.
An operator between Banach spaces is $ckw$-to-norm continuous if
and only if it is {\it completely continuous,} i.e., it
takes weakly null sequences into norm null
sequences \cite[Proposition~3.2]{GGkw}.

    If a Banach space contains a copy of $\ell_1$, then its unit
ball will not be $ckw$-precompact. This is the case, for instance,
of $L_1[0,1]$.
Moreover, it is proved in \cite[Theorem~3.7]{GGkw} that a bounded
subset $K$ of $L_\infty [0,1] = L_1[0,1]^*$,
is an ($L$)-set if and only if it is relatively compact as a subset
of $L_1[0,1]$, when we consider $L_\infty
[0,1]$ embedded into $L_1[0,1]$ by means of the identity map.

    We first give a result whose proof uses standard techniques.
To each polynomial $P\in\pkXY$, we associate an operator
$$
T_P:X \lra {\cal L}(^{k-1}X,Y)
$$
given by
$$
T_P(x)(x_1,\ldots,x_{k-1}):=\hat{P}(x,x_1,\ldots,x_{k-1}),
$$
for $x, x_1,\ldots, x_{k-1}\in X$.

\begin{Prop}
\label{tauunif}
Let $\tau$ be a vector topology on a Banach space $X$, and $P\in\pkXY$.
Then $P$ is uniformly $\tau$-continuous on bounded sets if and only if
so is $T_P$.
\end{Prop}

\Proof
Suppose $P$ is uniformly $\tau$-continuous on $B_X$.
Given $\epsilon>0$, we can find a balanced $\tau$ zero neighbourhood
$U$ in $X$ so that $\| Px-Py\|<\epsilon$ whenever $x,y\in B_X$
satisfy $x-y\in U$.

Assume $x,y$ satisfy the above conditions, and let $z_2,\ldots,z_k
\in B_X$. By the polarization formula,
{\small
\be
\lefteqn{\left( T_P(x)-T_P(y)\right) (z_2,\ldots,z_k) =}\\
& & \hat{P}(x,z_2,\ldots,z_k) - \hat{P}(y,z_2,\ldots,z_k) =\\
& & \frac{k^k}{2^kk!} \sum_{\epsilon_j=\pm 1} \epsilon_1\cdots\epsilon_k
\left[ P\left(\frac{\epsilon_1x+\epsilon_2z_2+\cdots +\epsilon_kz_k}{k}
    \right) -
   P\left(\frac{\epsilon_1y+\epsilon_2z_2+\cdots +\epsilon_kz_k}{k}
    \right) \right] \, .
\ee
}
Easily, we conclude that
$$
\| T_P(x)-T_P(y)\| \leq \epsilon\,\frac{k^k}{k!}\; ,
$$
and $T_P$ is uniformly $\tau$-continuous on $B_X$.

    Conversely, let $T_P$ be uniformly $\tau$-continuous on $B_X$.
For $0<\epsilon <1$, there is a $\tau$ zero neighbourhood $U\subset X$
so that $\| T_P(x)-T_P(y)\|<\eps$ whenever $x,y\in B_X$ satisfy $x-y\in U$.
For such $x,y$, we have
{\small
\be
\lefteqn{\| Px-Py\| \leq}\\
& &\| \hat{P}(x,\ldots,x) - \hat{P}(x,y,x,\ldots,x)\|
    + \|\hat{P}(x,y,x,\ldots,x) - \hat{P}(x,y,y,x,\ldots,x)\| +\cdots\\
& & \mbox{} + \|\hat{P}(x,y,\ldots,y) - \hat{P}(y,\ldots,y)\| =\\
& & \|\left( T_P(x)-T_P(y)\right) (x,\ldots,x)\| +
    \| \left( T_P(x)-T_P(y)\right) (x,y,x,\ldots,x)\| +\cdots <\\
& & k\epsilon \, ,
\ee
}
and $P$ is uniformly $\tau$-continuous on $B_X$, which completes the
proof.
\finesp

   Since $T_P$ is linear, $T_P$ is uniformly $\tau$-continuous on bounded
sets if and only if it is $\tau$-continuous on bounded sets.
Therefore, we easily obtain:

\begin{Cor}
\label{ckwunif}
A polynomial $P\in\pkXY$ is uniformly $ckw$-continuous on bounded sets
if and only if the associated operator $T_P$ is completely continuous.
\end{Cor}

    We are now ready to construct our polynomial, which is a
modification of an example given in \cite{ACL}.
    To this end, we divide
the interval $(0,1)$ into subintervals
$$
I_j:= \left( \frac{1}{2^j}\, ,\,\frac{1}{2^{j-1}}\right)
\hspace{2em} (j\in\N ),
$$
and denote by $\chi_j$ the characteristic function of $I_j$.
Let $(r_n)$ be the Rademacher functions on $[0,1]$, given by
$$
r_n(t):= \mbox{sign} \sin 2^n\pi t \hspace{2em} \left( t\in [0,1]\right)\, .
$$
Easily, we obtain:
\be
\langle r_n,\chi_j\rangle = \left\{ \begin{array}{cl}
              0 &  \hspace{1.5em}n>j \\
              -2^{-j} & \hspace{1.5em} n=j\\
              2^{-j} & \hspace{1.5em} n<j\, .
              \end{array} \right.
\ee
Note that for each $f\in L_1[0,1]$ we have
\begin{eqnarray}
\label{sumaf}
\sum_{j=1}^\infty
\left| \langle f,\chi_j\rangle\right| \leq \| f\|\, .
\end{eqnarray}

For $f, g, h\in L_1[0,1]$, define
$$
A(f,g,h) := \sum_{j=1}^\infty \langle f,\chi_j\rangle \langle g,\chi_j
\rangle\langle h,r_j\rangle\, .
$$
Since
$$
|A(f,g,h)|\leq \| g\|\cdot\| h\|\cdot \sum_{j=1}^\infty \left|\langle
f,\chi_j\rangle\right| \leq \| f\|\cdot\| g\|\cdot\| h\| \, ,
$$
we have that
$A$ is a $3$-linear continuous form on $L_1[0,1]$. Then, the function
$$
P(f):= A(f,f,f)\hspace{2em}\left( f\in L_1[0,1]\right)
$$
is a $3$-homogeneous continuous polynomial on $L_1[0,1]$.

\begin{Prop}
The polynomial $P$ is $ckw$-continuous on bounded sets of $L_1[0,1]$,
but is not uniformly $ckw$-continuous on bounded sets.
\end{Prop}

\Proof
We first show that the associated operator $T_P$ is not completely
continuous. By Corollary~\ref{ckwunif},
this will imply that $P$ is not uniformly $ckw$-continuous on bounded
sets.
Since
$$
\hat{P} (f,g,h)= \frac{1}{3}\left[ A(f,g,h)+A(h,f,g)+A(g,h,f)\right]\, ,
$$
the operator $T_P:L_1[0,1]\ra {\cal L}\left( ^2L_1[0,1]\right)$
is given by
$$
T_P(f)=\frac{1}{3}\left[
    A(\,\cdot\, ,\,\cdot\, ,f)
+ A(f,\,\cdot\, ,\,\cdot\, ) +
    A(\,\cdot\, ,f,\,\cdot\, )
    \right]\, .
$$
Then
\be
3\, T_P(r_n) & = & \sum_{j=1}^\infty \left[ \langle\,\cdot\, ,
\chi_j\rangle^2\langle r_n,r_j\rangle + 2\langle\,\cdot\, ,
\chi_j\rangle\langle\,\cdot\, ,r_j\rangle\langle r_n,\chi_j\rangle
\right] \\
& = & \langle\,\cdot\, ,\chi_n\rangle^2 -2\cdot
2^{-n}\langle\,\cdot\, ,\chi_n\rangle\langle\,\cdot\, ,r_n\rangle +
2\sum_{j=n+1}^\infty 2^{-j}\langle\,\cdot\, ,\chi_j\rangle
\langle\,\cdot\, ,r_j\rangle\; .
\ee
Since
$$
\left\| \langle\,\cdot\, ,\chi_n\rangle^2\right\| = 1\, ,
$$
and
$$
\left\| -2\cdot 2^{-n}\langle\,\cdot\, ,\chi_n\rangle\langle\,\cdot\,
,r_n\rangle + 2\sum_{j=n+1}^\infty 2^{-j}\langle\,\cdot\, ,\chi_j
\rangle\langle\,\cdot\, ,r_j\rangle\right\| \leq 4\cdot 2^{-n}\, ,
$$
we obtain that $\| T_P(r_n)\|$ does not converge to zero.
Since $(r_n)$ is a weakly null sequence in $L_1[0,1]$, this shows that
$T_P$ is not completely continuous.

Let us now prove that $P$ is $ckw$-continuous on bounded sets, in other
words, given $f_0\in L_1[0,1]$, $\eps >0$ and $r\geq \| f_0\|$,
there exists an ($L$)-set $K\subset L_1[0,1]^*$ so that
$|P(f)-P(f_0)|<\eps$ whenever $f\in L_1[0,1]$ satisfies $\| f\| \leq
r$ and $f-f_0\in\hspace{.3em} ^\circ\! K$, where
$$
^\circ\!K:= \left\{ g\in L_1[0,1]:|\langle g,h\rangle|\leq 1
\mbox{ for all
$h\in K$}\right\} \, .
$$
Indeed, choose $\delta >0$ with
\begin{eqnarray}
\label{delta}
5r^2\delta <\eps\, .
\end{eqnarray}
Thanks to (\ref{sumaf}), we can find $j_0\in\N$ so that
\begin{eqnarray}
\label{jota}
\left|\langle f_0,\chi_j\rangle \right| <\delta \hspace{1em}
\mbox{ for every $j>j_0$}.
\end{eqnarray}
The set
$$
K:= \left\{ \delta^{-1}\chi_j:j\in {\bf N}\right\} \bigcup
\left\{ \delta^{-1}r_j: 1\leq j\leq j_0\right\}
$$
is an ($L$)-set in $L_1[0,1]^*$.
Moreover, for $f\in L_1[0,1]$, we can write
\be
\lefteqn{P(f)-P(f_0)=}\\
& & A(f,f-f_0,f) + A(f,f_0,f-f_0) + A(f-f_0,f_0,f_0) =\\
& & \sum_{j=1}^\infty \langle f,\chi_j\rangle\langle f-f_0,\chi_j\rangle
    \langle f,r_j\rangle
   +\sum_{j=1}^{j_0} \langle f,\chi_j\rangle\langle f_0,\chi_j\rangle
    \langle f-f_0,r_j\rangle +\\
& & \sum_{j=j_0+1}^\infty \langle f,\chi_j\rangle
                       \langle f_0,\chi_j\rangle
    \langle f-f_0,r_j\rangle +
    \sum_{j=1}^\infty\langle f-f_0,\chi_j\rangle\langle f_0,\chi_j\rangle
    \langle f_0,r_j\rangle \; .
\ee
If $f-f_0\in \hspace{.3em}^\circ\! K$, we have
\be
\left|\langle f-f_0,r_j\rangle \right| \leq \delta
\hspace{1em}\mbox{ for } 1\leq j\leq
         j_0, \mbox{ and}\\
\left|\langle f-f_0,\chi_j\rangle \right| \leq \delta
\hspace{1em}\mbox{ for all } j\in\N\, .
\ee
Using these two inequalities along with (\ref{sumaf}),
(\ref{delta}) and (\ref{jota}), we get
\be
\lefteqn{|P(f)-P(f_0)| \leq }\\
& & \| f\|\cdot\delta\cdot\| f\| + \| f\|\cdot\| f_0\|\cdot\delta +
\| f\|\cdot \delta\cdot\| f-f_0\| + \delta\cdot \| f_0\|^2 \leq \\
& & 5r^2\delta < \eps \, ,
\ee
and the proof is complete.
\finesp

When we consider the norm topology, we do have that if a polynomial
is norm continuous (on bounded sets), it is also uniformly continuous
on bounded sets, and however the unit ball is not norm precompact in
infinite dimensional Banach spaces.
This is due to the fact that the balls centered at the origin constitute
a base of zero neighbourhoods, and a polynomial on a locally convex
space is continuous if and only if it is uniformly continuous on
some zero neighbourhood \cite{Dinbook}.

\section{Factorization results}
\label{factres}

In this part, we
apply Theorem~\ref{klinear} to prove that a holomorphic mapping $f$
is weakly uniformly continuous on bounded sets if and only if
it may be written in the form $f=g\circ S$, where $S$ is a compact operator,
and $g$ a holomorphic mapping.

For Banach spaces $X_1,\ldots,X_k$ and $Y$, we use
$\LXkY$ to represent the space of $k$-linear
(continuous) mappings from $X_1\times\cdots\times X_k$ into $Y$.
The space of compact operators from $X$ into $Y$ is denoted by
${\cal C}o(X;Y)$.

We denote by $\LwbXkY$ the space of $k$-linear mappings which are
weakly continuous on bounded sets, in the sense of Definition~\ref{bwcont}.
We define the spaces $\LwbkXY$ and $\pwbkXY$ in an analogous way.

For each $1\leq i\leq k$, consider the mapping
$$
\xi_i:\LXkY \lra {\cal L}(X_1,\ldots,X_i;{\cal L}(X_{i+1},\ldots ,X_k;Y))
$$
taking $A$ into $\overline{A}$ given by
$$
\overline{A}(x_1,\ldots,x_i)(x_{i+1},\ldots,x_k):= A(x_1,\ldots,x_k)
$$
for $x_j\in X_j$ $(1\leq j\leq k)$.
It is well known that $\xi_i$ is a linear surjective isometry.

\begin{Prop}
\label{eqsp}
Given Banach spaces $X_1,\ldots,X_k,Y$,
the isometry $\xi_i$ maps the space $\LwbXkY$ onto the space
$$
{\cal L}_{wb}\left( X_1,\ldots,X_i;{\cal L}_{wb}(X_{i+1},\ldots,X_k;Y)\right)\,
.
$$
\end{Prop}

\Proof
If $A$ is weakly continuous on bounded sets, it is clear that the range
of $\overline{A}$ lies in ${\cal L}_{wb}(X_{i+1},\ldots ,X_k;Y)$.
Now, let bounded, weak Cauchy nets $\left( x^j_\alpha\right) \subset X_j$
be given, for $1\leq j\leq i$, at least one of which being weakly null.
Suppose we have
$$
\left\| \overline{A} \left( x^1_\alpha,\ldots,x^i_\alpha\right)\right\| >
\delta
$$
for some $\delta >0$. Then we can find nets $\left( x^j_\alpha \right)
\subset B_{X_j}$ $(i+1\leq j\leq k)$ that can be assumed to be weak Cauchy,
so that
$$
\left\| A\left( x^1_\alpha,\ldots,x^k_\alpha\right)\right\| =
\left\| \overline{A}\left( x^1_\alpha,\ldots,x^i_\alpha \right)
\left( x^{i+1}_\alpha,\ldots,x^k_\alpha \right)\right\| >\delta\, ,
$$
a contradiction (Theorem \ref{klinear}). Applying
again Theorem~\ref{klinear},
we obtain that $\overline{A}$ is weakly continuous on bounded sets.

Conversely, if
$$
\overline{A} \in {\cal L}_{wb}\left(
X_1,\ldots,X_i;{\cal L}_{wb}(X_{i+1},\ldots,
X_k;Y)\right)
$$
and, for each $1\leq j\leq k$, a bounded net $\left( x^j_\alpha \right)
\subset X_j$ is given which converges weakly to $x^j$,
we easily obtain that $A\in\LwbXkY$ by writing
\be
\lefteqn{\left\| A\left( x^1_\alpha,\ldots,x^k_\alpha\right)
- A(x^1,\ldots,x^k)\right\| }\\
& \leq & \left\| \overline{A}\left( x^1_\alpha,\ldots,x^i_\alpha\right)
\left( x^{i+1}_\alpha,\ldots,x^k_\alpha\right)
-
            \overline{A}(x^1,\ldots,x^i)\left(
x^{i+1}_\alpha,\ldots,x^k_\alpha\right)\right\|
\\
& + & \left\|\overline{A}(x^1,\ldots,x^i)\left(
x^{i+1}_\alpha,\ldots,x^k_\alpha\right)
-
           \overline{A}(x^1,\ldots,x^i)(x^{i+1},\ldots,x^k)\right\| \, .
\ee
This completes the proof.
\finesp

Our next result shows that $P$ factorizes if and only if so does $\hat{P}$.

\begin{Prop}
\label{ideal}
Let ${\cal U}$ be an operator ideal, and $P\in\pkXY$, for Banach spaces
$X$, $Y$. The following
assertions are equivalent:

{\rm (a)} there are a Banach space $Z$, an operator $S\in {\cal U}(X;Z)$, and
a polynomial $Q\in\pk Z;Y)$ so that $P=Q\circ S$;

{\rm (b)} there are Banach spaces $Z_i$, operators $S_i\in {\cal U}(X;Z_i)$
$(1\leq i \leq n)$, and a $k$-linear mapping $B\in\LZkY$ so that
$\hat{P} =B\circ (S_1,\ldots ,S_k)$.

Moreover, if {\rm (b)} is satisfied, then we can choose $S$ and $Q$ in
{\rm (a)} so that $\| S\| = \max \| S_i\|$, and $\| Q\|=\| B\|$.
\end{Prop}

\Proof
(a) $\Ra$ (b). Take $Z_i=Z$, $S_i=S$ $(1\leq i\leq k)$, and $B=\hat{Q}$.

(b) $\Ra$ (a). Take $Z=Z_1\times\cdots\times Z_k$; $Sx:=(S_1x,\ldots ,S_kx)$
for all $x\in X$, and
$$
Q\left( (z_1,\ldots,z_k)\right) : = B(z_1,\ldots,z_k) \, .
$$
To see that $Q$ is a polynomial, note that its associated symmetric $k$-linear
mapping is
$$
\hat{Q}\left( \left( z_1^1,\ldots,z_k^1\right) ,\ldots ,\left( z_1^k,\ldots ,
z_k^k\right) \right) := \frac{1}{k!} \sum B\left( z_1^{i_1},\ldots ,
z_k^{i_k}\right) \, ,
$$
where the sum is taken over all permutations $(i_1,\ldots,i_k)$ of
the numbers $(1,\ldots,k)$.

Endowing $Z$ with the supremum norm, we obtain $\| S\| = \max\| S_i\|$,
and $\| Q\| =\| B\|$.
\finesp

    In the proof of the next theorem, we shall use the well known fact
that for every compact operator $T:X\ra Y$ between Banach spaces, we
can find a space $Z$ and compact operators $S:X\ra Z$ and $R:Z\ra Y$
so that $T=R\circ S$
\cite[Theorem~17.1.4]{J}.
We are indebted to H. Junek who pointed out that, as shown in
\cite[Lemma~1.2]{H}, we can renorm $Z$ so that
$\|R\|\!\cdot\!\|S\| = \|T\|$.
As usual, we denote by $c_0(X)$ the Banach space of all null sequences
in $X$, endowed with the supremum norm.

\begin{Th}
\label{facklin}
Given Banach spaces $X_1,\ldots,X_k,Y$ and a number
$\epsilon > 0,$ for each
$A\in\LwbXkY$ there are Banach spaces $Z_1,\ldots,Z_k$, operators
$S_i\in{\cal C}o(X_i;Z_i)$ and a mapping $B\in\LwbZkY$ so that
$$
A(x_1,\ldots,x_k) = B(S_1x_1,\ldots,S_kx_k) \hspace{2em} (x_i\in X_i)\, ,
$$
and
$\| B\|\cdot\| S_1\|\cdot\ldots\cdot\|S_k\| \leq (1+\epsilon )
\| A\|$.
\end{Th}

\Proof
For $k=1$, we have the above mentioned result for linear operators.
Assume the theorem is true for the $(k-1)$-linear
mappings $(k>1)$, and take $\lambda := (1+\epsilon )^{1/k}.$
Given $A\in\LwbXkY$, Proposition~\ref{eqsp} provides an associated
$$
\overline{A}\in{\cal L}_{wb}(X_2,\ldots,X_k;{\cal C}o(X_1;Y)) \, .
$$
By induction, we can write
$\overline{A} = \overline{D}\circ (S_2,\ldots,S_k)$, with $S_i\in
{\cal C}o(X_i;Z_i)$ $(2\leq i\leq k)$, $\overline{D}\in
{\cal L}_{wb}(Z_2,\ldots,Z_k;{\cal C}o(X_1;Y))$, and
$$
\|\overline{D}\|\cdot\|S_2\|\cdot\ldots\cdot\|S_k\| \leq
\lambda^{k-1}\|\overline{A}\|\, .
$$
By Proposition~\ref{eqsp}, we associate to
$\overline{D}$ an operator
$$
D\in {\cal C}o\left( Z_2;{\cal C}o\left( X_1;{\cal L}_{wb}(Z_3,\ldots,Z_k;
Y)\right) \right) \, .
$$
Since $D$ is compact, there is a sequence $(D_n)\subset
{\cal C}o\left( X_1;{\cal L}_{wb}(Z_3,\ldots,Z_k;
Y)\right)$ with $\| D_n\|\ra 0$ so that $D\left( B_{Z_2}\right)$ is
contained in the absolutely convex, closed hull of $\{ D_n\}$,
and
$$
\lambda^{-1} \leq \frac{\|D\|}{\sup \|D_n\|} \; .
$$
Define
$$
T:X_1\lra c_0 \left(
{\cal L}_{wb}(Z_3,\ldots,Z_k;
Y)\right)
$$
by  $Tx_1:= (D_nx_1)_{n=1}^\infty$.
Clearly, $T$ is compact and so we can find a space $Z_1$ and compact
operators $S_1:X_1\ra Z_1$ and $R:Z_1\ra
c_0 \left(
{\cal L}_{wb}(Z_3,\ldots,Z_k;
Y)\right)$ with
$S_1(X_1)$ dense in $Z_1$, such that $T=R\circ S_1$ and
$\|R\|\cdot\|S_1\| =    \|T\|$.
Define a linear mapping
$$
U:T(X_1)\lra {\cal L}_{wb}(Z_2,\ldots,Z_k;Y)
$$
by
$$
U(Tx_1)(z_2,\ldots,z_k):= (Dz_2)(x_1)(z_3,\ldots,z_k)  \, ,
$$
for $x_1\in X_1$, $z_2\in Z_2,\ldots,z_k\in Z_k$.
Clearly, $U$ is well defined.
If $\| z_2\| = 1$, we have $Dz_2 = \sum_{n=1}^\infty \lambda_nD_n$,
with $\sum_{n=1}^\infty |\lambda_n| \leq 1$.
Then,
\be
\left\| U(Tx_1)(z_2,\ldots,z_k)\right\|
& = & \left\| \sum_{n=1}^\infty
       \lambda_n(D_nx_1)(z_3,\ldots,z_k)\right\| \\
& \leq & \left( \sum_{n=1}^\infty |\lambda_n|\right)
         \left( \sup_n \|D_nx_1\|
              \right) \cdot \| z_3\|\cdot\ldots\cdot\| z_k\| \\
& \leq & \| Tx_1\|\cdot \| z_2\|\cdot\ldots\cdot\| z_k\|\, .
\ee
Therefore, $U$ is continuous and admits an extension $V$ to the
closure of $T(X_1)$, with $\|V\| =\|U\|\leq 1$.
Since $R$ is compact, the operator
$$
V\circ R: Z_1\lra {\cal L}_{wb}(Z_2,\dots,Z_k;Y)
$$
is compact. Let
$B\in {\cal L}_{wb}(Z_1,\dots,Z_k;Y)$ be the $k$-linear mapping
associated to $V\circ R$ by Proposition~\ref{eqsp}.
We have
$$
\| B\| = \|V\circ R\| \leq \|R\| \, .
$$
For $x_1\in X_1,\ldots,x_k\in X_k$, we obtain
\be
B(S_1x_1,\ldots,S_kx_k) & = & \left( V\circ R (S_1x_1)\right)
         (S_2x_2,\ldots,S_kx_k) \\
& = & U(Tx_1) (S_2x_2,\ldots,S_kx_k) \\
& = & (DS_2x_2)(x_1)(S_3x_3,\ldots,S_kx_k) \\
& = & \overline{D} (S_2x_2,\ldots,S_kx_k)(x_1) \\
& = & \overline{A}(x_2,\ldots,x_k)(x_1) \\
& = & A(x_1,\ldots ,x_k) \, .
\ee
Moreover, since $\|\overline{A}\|=\|A\|$ and
$$
\|\overline{D}\| =\|D\| \geq \lambda^{-1} \sup\|D_n\|
=\lambda^{-1}\|T\|=\lambda^{-1} \|R\|
\cdot\|S_1\| \geq \lambda^{-1} \|B\|\cdot\|S_1\|\, ,
$$
we get
$$
\|B\|\cdot\|S_1\|\cdot\|S_2\|\cdot\ldots\cdot\|S_k\|
\leq \lambda \|\overline{D}\|\cdot\|S_2\|\cdot\ldots\cdot\|S_k\|
\leq \lambda^{k} \|A\|\, ,
$$
and the proof is complete.
\fin

\begin{Cor}
\label{facpol}
Given Banach spaces $X, Y$ and a polynomial
$P\in\pkXY$, we have that
$P\in\pwbkXY$ if and only if there are a Banach space $Z$, an operator
$S\in{\cal C}o(X;Z)$, and a polynomial $Q\in\pwbk Z;Y)$
such that $P=Q\circ S$.
Moreover, given $\eps >0$,
we can obtain $\| Q\|\cdot \| S\|^k \leq (1+\epsilon)
\|\hat{P}\|$.
\end{Cor}

\Proof
Suppose $P\in\pwbkXY$. Then
$\hat{P}\in\LwbkXY$. Given $\eps >0$,
by Theorem~\ref{facklin}, we can write $\hat{P}=B\circ (S_1,\ldots,S_k)$,
with $S_i\in{\cal C}o(X;Z_i)$, $B\in\LwbZkY$, and
$$
\| B\|\cdot\|S_1\|\cdot\ldots\cdot\|S_k\| \leq (1+\eps )
\|\hat{P}\| \, .
$$
We can assume that $\|S_1\| = \cdots = \|S_k\|$.
By Proposition~\ref{ideal}, we have $P=Q\circ S$, with
$\|Q\| \cdot \|S\|^k \leq (1+\eps ) \|\hat{P}\|$.
Easily, $Q$ is weakly continuous on bounded sets.
The converse is clear.
\finesp

      For complex Banach spaces $X$, $Y$, let $\HXY$ denote the
space of all holomorphic mappings from $X$ into $Y$, and
$\HwbuXY$ the subspace
of all mappings in $\HXY$ whose restrictions to
bounded subsets are weakly uniformly continuous.
We obtain:

\begin{Th}
Let $X$, $Y$ be complex Banach spaces, and
$f\in\HXY$. Then $f\in\HwbuXY$ if and only if there are a space $Z$,
an operator $S\in{\cal C}o(X;Z)$, and a mapping $g\in{\cal H}_{wbu}(Z;Y)$
so that $f=g\circ S$.
\end{Th}

\Proof
Let $f=\sum_{k=1}^\infty P_k$ be the Taylor series expansion of $f$
at the origin, and suppose $f\in\HwbuXY$.
By the Cauchy-Hadamard formula, we have $\lim\|P_k\|^{1/k} =0$.
By Corollary~\ref{facpol}, there are spaces $Z_k$,
operators $S_k\in{\cal C}o(X;Z_k)$
and polynomials $Q_k\in\pwbk Z_k;Y)$ such that $P_k=Q_k\circ S_k$, with
$$
\|Q_k\|\cdot\|S_k\|^k \leq 2 \|
\hat{P}_k\| \leq 2\frac{k^k}{k!}\,\|P_k\|
$$
(the last inequality is well known and
may be seen in \cite[Theorem~2.2]{Mu}).
Then, using the Stirling formula,
$$
\lim_k \|Q_k\|^{1/k}\cdot \|S_k\| \leq \lim_k
\frac{2^{1/k}e}{(2\pi k)^{1/2k}}\cdot
\|P_k\|^{1/k}   =0 \, .
$$
We can assume therefore that $\|S_k\|\ra 0$ and $\|Q_k\|^{1/k}\ra 0$.
Define
$$
S: X \lra Z:= c_0(Z_k)
$$
by $Sx= (S_kx)_k$.
Clearly, $S$ is compact.
Denoting
$$
\pi_k: (y_i)\in Z \longmapsto y_k\in Z_k \, ,
$$
we define $g:Z\ra Y$ by $g(y):= \sum_{k=1}^\infty Q_k\circ \pi_k(y)$.
Since $\lim\|Q_k\circ\pi_k\|^{1/k} = \lim\|Q_k\|^{1/k} =0$, we have
that $g$ is a holomorphic mapping,
bounded on bounded sets. Moreover, $Q_k\circ \pi_k$ is
weakly continuous on bounded sets, for all $k$. Therefore, $g\in{\cal H}_{wbu}
(Z;Y)$. The converse is clear.
\finesp

We recall that it remains unknown whether a holomorphic mapping
between Banach spaces which is weakly continuous on bounded sets
is or not automatically weakly uniformly continuous
on bounded sets \cite{AHV}, i.e., whether or not a holomorphic
function on a Banach space can satisfy the conditions (a), (b) and
(c) given after the proof of Theorem~\ref{equivAHV}.\espvv

We thank Professor S. Dineen for pointing out a gap in a first version of
Proposition~\ref{eqsp}. Trying to fill in this gap led us to find
Theorem~\ref{klinear}.

Most of this work was done while the second named author was visiting the
Mathematics Department of University College Dublin, whom
he wishes to thank for their hospitality.
He also thanks DGICYT (Spain), for supporting his stay.

\espvv
\noin
{\small Manuel Gonz\'alez \\
      Departamento de Matem\'aticas \\
      Facultad de Ciencias\\
      Universidad de Cantabria \\ 39071 Santander (Spain)\\
      e-mail: gonzalem@ccaix3.unican.es\espvv \\
      Joaqu\'\i n  M. Guti\'errez\\
      Departamento de Matem\'atica Aplicada\\
      ETS de Ingenieros Industriales \\
      Universidad Polit\'ecnica de Madrid\\
      C. Jos\'e Guti\'errez Abascal 2 \\
      28006 Madrid (Spain)\\
      e-mail: c0550003@ccupm.upm.es}

\end{document}